\documentclass{ article}

\begin{document}
\vspace*{.5cm}
\begin{center}
{\Large{\bf   Anti-invariant Riemannian maps from almost Hermitian manifolds}}\\
\vspace{.5cm}
 { Bayram \d{S}ahin} \\
\end{center}

\vspace{.5cm}
\begin{center}
{\it Inonu University, Department of Mathematics, 44280,
Malatya-Turkey. E-mail:bayram.sahin@inonu.edu.tr}
\end{center}
\vspace{.5cm} \noindent {\bf Abstract.} {\small As a generalization of anti-invariant Riemannian submersions, we introduce
anti-invariant Riemannian maps  from almost Hermitian manifolds to
Riemannian manifolds. We give examples and investigate the geometry
of foliations which are arisen from the definition of an
anti-Riemannian map. Then we give a decomposition theorem for the
source manifold of such maps. We also find necessary and sufficient
conditions for anti-invariant Riemannian maps to be totally geodesic and show that every pluriharmonic Lagrangian Riemannian map, which is a special anti-invariant Riemannian map, from a K\"{a}hler manifold to a Riemannian manifold is totally geodesic.}\\

 \noindent\small{{\bf 2000 Mathematics Subject Classification:}
53C15, 58C25, 53C43.\\

}

\section*{\bf 1.~ Introduction}

\setcounter{equation}{0}
\renewcommand{\theequation}{1.\arabic{equation}}
\markboth{Anti-invariant maps}{Anti-invariant
maps}{\thispagestyle{plain}} Smooth maps between Riemannian
manifolds are useful for comparing geometric structures between two
manifolds. Isometric immersions (Riemannian submanifolds) are basic
such maps between Riemannian manifolds and they are characterized by
their Riemannian metrics and Jacobian matrices. More precisely, a
smooth map $F:(M_1, g_1)\longrightarrow (M_2, g_2)$ between
Riemannian manifolds $(M_1, g_1)$ and $(M_2, g_2)$ is called an
isometric immersion if $F_*$ is injective and
\begin{equation}
g_2(F_*X, F_*Y)=g_1(X, Y)\label{eq:1.1}
\end{equation} for $X, Y $
vector fields tangent to $M_1$, here $F_*$ denotes the derivative
map.\\

Let $\bar{M}$ be a K\"{a}hler manifold with complex structure $J$
and $M$ a Riemannian mani-
fold isometrically immersed in $\bar{M}$.
We note that  submanifolds of a K\"{a}hler manifold are determined
by the behavior of the tangent bundle of the submanifold under the
action of the complex structure of the ambient manifold. A
submanifold $M$ is called holomorphic (complex) if $J (T_p M)\subset
T_p M$, for every $p \in M$, where $T_p M$  denotes the tangent
space to $M$ at the point $p $. $M$ is called totally real if $J(T_p
M) \subset T_p M^{\perp}$ for every $p \in M,$ where $T_p M^{\perp}$
denotes the normal space to $M$ at the point $p$.\\

A smooth map $F:(M_1, g_1)\longrightarrow (M_2, g_2)$ is called a
Riemannian submersion if $F_*$ is onto and it satisfies the equation
(\ref{eq:1.1}) for vector fields tangent to the horizontal space $(kerF_*)^\perp$. Riemannian submersions between Riemannian manifolds
were studied by O'Neill \cite{O'Neill} and Gray \cite{Gray}, see
also \cite{Falcitelli-Ianus-Pastore}.  Later such submersions were
considered between manifolds equipped with differentiable structures. As an
analogue of holomorphic submanifolds, Watson defined almost
Hermitian submersions between almost Hermitian manifolds and he
showed that the base manifold and each fiber have the same kind of
structure as the total space, in most cases \cite{Watson}. We note
that almost Hermitian submersions have been extended to the almost
contact manifolds \cite{Domingo}, locally conformal K\"{a}hler
manifolds \cite{Marrero-Rocha} and quaternion K\"{a}hler manifolds
\cite{Ianus-Mazzocco-Vilcu}.\\

Let $M$ be a complex $m-$dimensional Hermitian manifold with
Hermitian metric $g_M$ and almost complex structure $J_M$ and $N$ a
complex $n-$dimensional Hermitian manifold with Hermitian metric
$g_N$ and almost complex structure $J_N$. A Riemannian submersion
$F:M\longrightarrow N$ is called an almost Hermitian submersion if
$F$ is an almost complex map, i.e., $F_*J_M=J_NF_*$. The main result
of this notion is that the vertical and horizontal distributions are
$J_M-$ invariant. On the other hand, Riemannian submersions from
almost Hermitian manifolds onto Riemannian manifolds have been
studied by many authors under the assumption that the vertical
spaces of such submersions are invariant with respect to the complex
structure.
 For instance, Escobales \cite{Escobales} studied Riemannian
submersions from complex projectives space onto Riemannian manifolds
under the assumption that the fibers are connected, complex, totally
geodesic submanifolds. One can see that this assumption implies that
the vertical distribution is invariant.\\

Recently, we introduced anti-invariant Riemannian submersions from
almost Hermitian manifolds onto Riemannian manifolds  and
investigated the geometry of such submersions \cite{Sahin1}. We now
recall the definition of anti-invariant Riemannian submersions. Let
$M$ be a complex $m-$ dimensional almost Hermitian manifold with
Hermitian metric $g_{_M}$ and almost complex structure $J$ and $N$
be a Riemannian manifold with Riemannian metric $g_{_N}$. If there exists a Riemannian submersion $F:M \longrightarrow N$
such that $ker F_*$ is anti-invariant with respect to $J$,i.e.,
$J(ker F_*)\subseteq (ker F_*)^{\perp}$, then we say that $F$ is an
anti-invariant Riemannian submersion.\\

 In 1992, Fischer  introduced
Riemannian maps between Riemannian manifolds in \cite{Fischer} as
a generalization of the notions of isometric immersions and
Riemannian submersions (For Riemannian maps and their applications
in spacetime geometry, see; \cite{Garcia-Rio-Kupeli}). Let
$F:(M_1, g_1)\longrightarrow (M_2, g_2)$ be a smooth map between
Riemannian manifolds such that $0<rank F<min\{ m, n\}$, where
$dimM_1=m$ and $dimM_2=n$. Then we denote the kernel space of
$F_*$ by $kerF_*$ and consider the orthogonal complementary space
$\mathcal{H}=(kerF_*)^\perp$ to $kerF_*$. Suppose that $rank F$ is
constant. Then the tangent bundle of $M_1$ has the following
decomposition

$$TM_1=kerF_* \oplus\mathcal{H}.$$

We denote the range of $F_*$ by $rangeF_*$ and consider the
orthogonal complementary space $(rangeF_*)^\perp$ to $rangeF_*$
in the tangent bundle $TM_2$ of $M_2$. Since $rankF<min\{ m, n\}$,
we always have $(rangeF_*)^\perp\neq \{0\}$. Thus the tangent bundle
$TM_2$ of $M_2$ has the following decomposition
$$TM_2=(rangeF_*)\oplus (rangeF_*)^\perp.$$

 Now, a smooth
map $F:(M^{^m}_1,g_1)\longrightarrow (M^{^n}_2, g_2)$ is called
Riemannian map at $p_1 \in M$ if the horizontal restriction
$F^{^h}_{*p_1}: (ker F_{*p_1})^\perp \longrightarrow (range
F_{*p_1})$  is a line-
ar isometry between the inner product spaces
$((ker F_{*p_1})^\perp, g_1(p_1)\mid_{(ker F_{*p_1})^\perp})$ and
$(range F_{*p_1}, g_2(p_2)\mid_{(range F_{*p_1})})$, $p_2=F(p_1)$.
Therefore Fischer stated in \cite{Fischer} that a Riemannian map is
a map which is as isometric as it can be. In other words, $F_*$
satisfies the equation (\ref{eq:1.1}) for $X, Y$ vector fields
tangent to $\mathcal{H}$. It follows that isometric immersions and
Riemannian submersions are particular Riemannian maps with
$kerF_*=\{ 0 \}$ and $(rangeF_*)^\perp=\{ 0 \}$. It is known that
a Riemannian map is a subimmersion\cite{Fischer}. It is also known that a map $F:M_1 \longrightarrow M_2$ is a subimmersion if and only if the rank of the linear map $F_{*p}:T_pM_1\longrightarrow T_{F(p)}M_2$ is constant for $p$ in each connected component of $M_1$ \cite{Marsden}.  \\

In this paper, as a generalization of anti-Riemannian submersions
from almost Hermitian manifolds to Riemannian manifolds, we
introduce anti-invariant Riemannian maps. The study is organized
as follows: In section~2, we present the basic information needed
for this paper. In section~3, we give definition of anti-invariant
Riemannian maps, provide examples and investigate the geometry of
leaves of the distributions. We then obtain a decomposition
theorem. Finally, in section 4, we obtain necessary and sufficient
conditions for such maps to be totally geodesic. We also show that
a Lagrangian Riemannian map with totally umbilical fibers implies
a totally geodesic foliation on $M_1$ and every pluriharmonic
Lagrangian map is totally geodesic.

\section*{2.~ Preliminaries}
\setcounter{equation}{0}
\renewcommand{\theequation}{2.\arabic{equation}}

In this section we recall some basic materials from
\cite{Baird-Wood}. Let $(M, g_{_M})$ be a Riemannian manifold and
$\mathcal{V}$ be a $q-$ dimensional distribution on $M.$ Denote its
orthogonal distribution $\mathcal{V}^{\perp}$ by $\mathcal{H}.$
Then, we have
\begin{equation}
TM=\mathcal{V}\oplus \mathcal{H}. \label{eq:2.1}
\end{equation}
$\mathcal{V}$ is called the vertical distribution and
$\mathcal{H}$ is called the horizontal distribution. We use the
same letters to denote the orthogonal projections onto these
distributions.

By the unsymmetrized second fundamental form of $\mathcal{V},$ we
mean the tensor field $A^{\mathcal{V}}$ defined by
\begin{equation}
A^{\mathcal{V}}_E F=\mathcal{H}(\nabla_{\mathcal{V}E} \mathcal{V}
F),\quad  E, F \in \Gamma(TM), \label{eq:2.2}
\end{equation}
where $\nabla$ is the Levi-Civita connection on $M.$ The
symmetrized second fundamental form $B^{\mathcal{V}}$ of
$\mathcal{V}$ is given by
\begin{equation}
B^{\mathcal{V}}(E, F)=\frac{1}{2} \{ A^{\mathcal{V}}_E F+
A^{\mathcal{V}}_F E\}=\frac{1}{2}\{
\mathcal{H}(\nabla_{\mathcal{V}E} \mathcal{V}
F)+\mathcal{H}(\nabla_{\mathcal{V}F} \mathcal{V}
E)\}\label{eq:2.3}
\end{equation}
for  any $ E, F \in \Gamma(TM).$
 The
integrability tensor of $\mathcal{V}$ is the tensor field
$I^{\mathcal{V}}$ given by
\begin{equation}
I^{\mathcal{V}}(E, F)=A^{\mathcal{V}}_E F-A^{\mathcal{V}}_F
E-\mathcal{H}([\mathcal{V}E, \mathcal{V}F]). \label{eq:2.4}
\end{equation}
Moreover, the mean curvature vector field of $\mathcal{V}$ is defined by
\begin{equation}
\mu^{\mathcal{V}}=\frac{1}{q}Trace
B^{\mathcal{V}}=\frac{1}{q}\sum^q_{i=1} \mathcal{H}(\nabla_{e_r}
e_r), \label{eq:2.5}
\end{equation}
where $\{ e_1,..., e_{q}\}$ is a local frame of $\mathcal{V}.$

 By reversing the roles of $\mathcal{V},$ $\mathcal{H},$
$B^{\mathcal{H}},$ $A^{\mathcal{H}}$ and $I^{\mathcal{H}}$ can be
defined similarly. For instance, $B^{\mathcal{H}}$ is defined by
\begin{equation}
B^{\mathcal{H}}(E, F)= \frac{1}{2}\{
\mathcal{V}(\nabla_{\mathcal{H}E} \mathcal{H}
F)+\mathcal{V}(\nabla_{\mathcal{H}F} \mathcal{H} E)\}
\label{eq:2.6}
\end{equation}
and, hence we have
\begin{equation}
\mu^{\mathcal{H}}=\frac{1}{m-q}Trace B^{\mathcal{H}}=\frac{1}{m-q}
\sum^{m-q}_{s=1} \mathcal{V}(\nabla_{E_s} E_s), \label{eq:2.7}
\end{equation}
where $E_1,..., E_{m-q}$ is a local frame of $\mathcal{H}$. A distribution $\mathcal{D}$ on $M$ is said to be  minimal if, for each $x \in M$, the mean curvature vector field vanishes. \\

Let $(M, g_{_M})$ and $(N, g_{_N})$ be Riemannian manifolds and
suppose that $\varphi: M\longrightarrow N$ is a smooth map
between them. Then the differential ${\varphi}_*$ of $\varphi$ can
be viewed a section of the bundle $Hom(TM,
\varphi^{-1}TN)\longrightarrow M,$ where $\varphi^{-1}TN$ is the
pullback bundle which has fibres
$(\varphi^{-1}TN)_p=T_{\varphi(p)} N, p \in M.$ $Hom(TM,
\varphi^{-1}TN)$ has a connection $\nabla$ induced from the
Levi-Civita connection $\nabla^M$ and the pullback connection.
Then the second fundamental form of $\varphi$ is given by
\begin{equation}
(\nabla {\varphi}_*)(X, Y)=\nabla^{\varphi}_X
{\varphi}_*(Y)-{\varphi}_*(\nabla^M_X Y) \label{eq:2.8}
\end{equation}
for $X, Y \in \Gamma(TM).$ It is known that the second fundamental
form is symmetric.\\

Finally, we recall the definition of almost Hermitian and
K\"{a}hler manifolds. Let ($\bar{M}, g$) be an almost Hermitian
manifold. This means \cite{Yano-Kon} that $\bar{M}$ admits a
tensor field $J$ of type (1,1) on $\bar{M}$ such that, $\forall X,
Y \in \Gamma(T\bar{M})$, we have
\begin{equation}
J^2=-I, \quad g(X, Y)=g(JX, JY) \label{eq:2.9}.
\end{equation}
 An almost Hermitian manifold $\bar{M}$ is called  K\"{a}hler manifold if
\begin{equation}
(\bar{\nabla}_XJ)Y=0,\forall X,Y \in \Gamma(T\bar{M}),
\label{eq:2.10}
\end{equation}
where $\bar{\nabla}$ is the Levi-Civita connection on $\bar{M}. $

\section*{\bf 3.~Anti-invariant Riemannian maps}
\setcounter{equation}{0}
\renewcommand{\theequation}{3.\arabic{equation}}
In this section, we define anti-invariant Riemannian maps, provide
examples and investigate the geometry of such maps. We first
present the following definition.\\

\noindent{\bf Definition~3.1.~}Let $F$ be a Riemannian map from an
almost Hermitian manifold $(M_1,g_1,J_1)$ to a Riemannian manifold
$(M_2,g_2)$. We say that $F$ is an anti-invariant Riemannian map
if the following condition is satisfied
$$J_1(ker F_*)\subset (ker F_*)^\perp.$$

We denote the orthogonal complementary subbundle  to $J_1(ker F_*)$
in $(ker F_*)^\perp$ by $\mu$. Then it is easy to see that $\mu$ is
invariant with respect to $J_1$.\\

Since $F$ is a subimmersion, it follows that the rank of $F$ is
constant on $M_1$, then the rank theorem for functions implies that
$ker F_*$ is an integrable subbundle of $TM_1$, (\cite{Marsden},
page:205).  Thus it follows from above definition that the leaves of
the distribution $ker F_*$ of a anti-invariant Riemannian map are totally real submanifolds of $M_1$.\\

\noindent{\bf Proposition~3.1.~}{\it Let $F$ be an anti-invariant
Riemannian map from an almost Hermitian manifold
$(M^{2m}_1,g_1,J_1)$ to a Riemannian manifold $(M^n_2,g_2)$ with
$dim(ker F_*)=r$. Then
$dim(\mu)=2(m-r)$.}\\

\noindent{\bf Proof.~} Since $dim(ker F_*)=r$ and $F$ is
anti-invariant Riemannian map, it follows that $dim(J_1(ker
F_*))=r$. Hence $dim(\mu)=2(m-r)$.\\

Let $F$ be an anti-invariant Riemannian map from an almost
Hermitian manifold $(M^{2m}_1,g_1,J_1)$ to a Riemannian manifold
$(M^n_2,g_2)$, then, as in the submanifold theory of Hermitian
manifolds, we say that $F$ is a Lagrangian Riemannian map from
$M_1$ to $M_2$ if $J_1(ker F_*)=(ker F_*)^\perp$. Hence  an
anti-invariant Riemannian map $F$ is a Lagrangian Riemannian map
if and only if $\frac{dim(M_1)}{2}=dim(ker F_*)$. It is also easy
to see that an anti-invariant Riemannian map is Lagrangian if and
only if $\mu=\{0\}$. Moreover, we have the
following.\\

\noindent{\bf Proposition~3.2.~}{\it Let $F$ be an anti-invariant
Riemannian map from an almost Hermitian manifold $(M_1,g_1,J_1)$
to a Riemannian manifold $(M_2,g_2)$ with $dim(ker F_*)=r$. Then
$F$ is a Lagrangian map if and only if
$\frac{dim(M_1)}{2}=dim(range F_*)$.}\\

We now give some examples for anti-invariant Riemannian maps.\\

\noindent{\bf Proposition~3.3.~}{\it Every proper Riemannian map $F$
($F$ is neither isometric immersion nor Riemannian submersion) from
an almost Hermitian manifold $(M^2_1,J,g_1)$ to a Riemannian
manifold $(M^n_2,g_2)$ is an anti-invariant Riemannian
map.}\\

\noindent{\bf Example~3.1.~}Every anti-invariant Riemannian
submersion from an almost Hermitian manifold to a Riemannian
manifold is an anti-invariant Riemannian map with $(range
F_*)^\perp=\{0\}$.\\

\noindent{\bf Example~3.2.~}Let $F$ be a map defined by
$$
\begin{array}{cccc}
  F: & R^4             & \longrightarrow & R^3\\
     & (x_1,x_2,x_3,x_4) &             & (\frac{x_1 -
     x_3}{\sqrt{2}},0,\frac{x_2+x_4}{\sqrt{2}}).
\end{array}
$$
Then it follows that
$$(ker F_*)=Span\{Z_1=\frac{\partial}{\partial x_1}+\frac{\partial}{\partial
x_3},\, Z_2=\frac{\partial}{\partial x_2}-\frac{\partial}{\partial
x_4}\}$$ and
$$(ker F_*)^\perp=Span\{Z_3=\frac{\partial}{\partial x_1}-\frac{\partial}{\partial
x_3},\, Z_4=\frac{\partial}{\partial x_2}+\frac{\partial}{\partial
x_4}\}.$$ By direct computations one can see that
$$F_*(Z_3)={\sqrt{2}}\frac{\partial}{\partial x_1},\,F_*(Z_4)={\sqrt{2}}\frac{\partial}{\partial
x_3}.$$ Thus $F$ is a Riemannian map with $(range
F_*)^\perp=Span\{\frac{\partial}{\partial x_2}\}$. Moreover it is
easy to see that $JZ_1=Z_4,\,JZ_2=-Z_3$, where $J$ is the canonical
complex structure of $R^4$ defined by \[J(x^1, \bar{x}^1, x^2,
\bar{x}^2)  =  (-\bar{x}^1, x^1, -\bar{x}^2, x^2).\] As a result,
$F$ is an
anti-invariant Riemannian map.\\

Let $F$ be anti-invariant Riemannian map from a K\"{a}hler
manifold $(M_1,g_1,J)$ to a Riemannian manifold $(M_2,g_2)$. Then
for $X \in \Gamma((ker F_*)^\perp)$, we have
\begin{equation}
JX=\mathcal{B}X+\mathcal{C}X, \label{eq:3.1}
\end{equation}
where $\mathcal{B}X \in \Gamma(ker F_*)$ and $\mathcal{C}X \in
\Gamma(\mu)$.\\

We now investigate the geometry of leaves of the distributions
$(ker F_*)$ and $(ker F_*)^\perp$.\\

\noindent{\bf Proposition~3.4.~}{\it Let $F$ be an anti-invariant
Riemannian map from a K\"{a}hler mani-
fold  $(M_1,g_1,J_1)$ to a
Riemannian manifold $(M_2,g_2)$. Then $(ker F_*)$ defines a totally
geodesic foliation on $M_1$ if and only if}
\begin{equation}
g_2((\nabla F_*)(X,\mathcal{B}Z),F_*(J_1Y))=g_2((\nabla
F_*)(J_1Y,X), F_*(\mathcal{C}Z))\label{eq:3.2}
\end{equation} {\it
for $X,Y \in \Gamma(ker F_*)$ and $Z\in \Gamma((ker F_*)^\perp)$.}

\noindent{\bf Proof.~} For $X, Y \in \Gamma(ker F_*)$ and $Z \in
\Gamma((ker F_*)^\perp)$, by using (\ref{eq:2.9}), (\ref{eq:2.10})
and (\ref{eq:3.1}) we have
$$g_1(\nabla^1_X Y,Z)=g_1(\nabla^1_XJ_1Y,
\mathcal{B}Z)+g_1(\nabla^1_XJ_1Y,\mathcal{C}Z),$$  where
$\nabla^1$ is the Levi-Civita connection on $M_1$. Since
$\mathcal{B}Z\in \Gamma(ker F_*)$, we get
$$g_1(\nabla^1_X Y,Z)=-g_1(J_1Y,
\nabla^1_X\mathcal{B}Z)+g_1(\nabla^1_XJ_1Y,\mathcal{C}Z).$$ Since
$F$ is a Riemannian map, using (\ref{eq:2.8}) we obtain
$$g_1(\nabla^1_X Y,Z)=g_2(F_*(J_1Y),
(\nabla F_*)(X,\mathcal{B}Z))-g_2((\nabla
F_*)(X,J_1Y),F_*(\mathcal{C}Z))$$ which proves the assertion.\\

For the distribution $(ker F_*)^\perp$, we have the following.\\

\noindent{\bf Proposition~3.5.~}{\it Let $F$ be an anti-invariant
Riemannian map from a K\"{a}hler manifold  $(M_1,g_1,J_1)$ to a
Riemannian manifold $(M_2,g_2)$. Then $(ker F_*)^\perp$ defines a
totally geodesic foliation on $M_1$ if and only if}
\begin{equation}g_2((\nabla
F_*)(Z_1,\mathcal{B}Z_2),F_*(J_1X))=-g_2(\nabla^{^F}_{Z_1}F_*(J_1X),F_*(\mathcal{C}Z_2))\label{eq:3.3}
\end{equation}
{\it for $Z_1,Z_2 \in \Gamma((ker F_*)^\perp)$ and $X\in
\Gamma(ker F_*)$.}\\

\noindent{\bf Proof.~} For $Z_1,Z_2 \in \Gamma((ker F_*)^\perp)$
and $X\in \Gamma(ker F_*)$, we have
$g_1(\nabla^1_{Z_1}Z_2,X)=-g_1(Z_2,\nabla^1_{Z_1}X)$. Using
(\ref{eq:2.9}), (\ref{eq:2.10}) and (\ref{eq:3.1}) we obtain
$$g_1(\nabla^1_{Z_1}Z_2,X)=-g_1(\nabla^1_{Z_1}J_1X,\mathcal{B}Z_2)-g_1(\nabla^1_{Z_1}J_1X,\mathcal{C}Z_2).$$
Hence, we get
$$g_1(\nabla^1_{Z_1}Z_2,X)=g_1(J_1X,\nabla^1_{Z_1}\mathcal{B}Z_2)-g_1(\nabla^1_{Z_1}J_1X,\mathcal{C}Z_2).$$
Then Riemannian map $F$ and (\ref{eq:2.8}) imply that
\begin{eqnarray}
g_1(\nabla^1_{Z_1}Z_2,X)&=&-g_2(F_*(J_1X),(\nabla
F_*)(Z_1,\mathcal{B}Z_2))+g_2((\nabla
F_*)(Z_1,J_1X),F_*(\mathcal{C}Z_2))\nonumber\\
&&-g_2(\nabla^F_{Z_1}F_*(J_1X), F_*(\mathcal{C}Z_2)).\nonumber
\end{eqnarray}
 From \cite{Sahin2}, we know that there are no components of
$(\nabla F_*)(Z_1,J_1X)$ in $(range F_*)$, hence we get the
result.\\

We now recall the de Rham theorem for Riemannian manifolds. Let $g$
be a Riemannian metric tensor on the manifold $M = M_1\times M_2$
and assume that the canonical foliations $L_{M_1}$ and $L_{M_2}$
intersect perpendicularly everywhere. Then $g$ is the metric tensor
of a usual product of Riemannian manifolds if and only if $L_{M_1}$
and $L_{M_2}$ are totally geodesic foliations \cite{Rham}. From
Proposition 3.4 and Proposition 3.5, we have the following
decomposition theorem.\\

\noindent{\bf Theorem~3.1.~}{\it Let $F$ be an anti-invariant
Riemannian map from a K\"{a}hler manifold $(M_1,g_1,J_1)$ to a
Riemannian manifold $(M_2,g_2)$. Then $M_1$ is a locally product
Riemannian manifold in the form $M_{(ker F_*)}\times M_{(ker
F_*)^\perp}$ if and only if (\ref{eq:3.2}) and (\ref{eq:3.3}) are
satisfied, where $M_{(ker F_*)}$ and $ M_{(ker
F_*)^\perp}$ are the integral manifolds of $kerF_*$ and $(ker F_*)^\perp$, respectively.}\\

\section*{\bf 4.Totally geodesic anti-invariant Riemannian maps}
\setcounter{equation}{0}
\renewcommand{\theequation}{4.\arabic{equation}}
In this section, we first obtain necessary and sufficient conditions for anti-invariant Riemannian maps from K\"{a}hler manifolds to Riemannian manifolds to be totally geodesic. Then we show that an anti-invariant Riemannian map with totally umbilical fibres has totally geodesic distribution $ker F_*$ and every pluriharmonic anti-invariant Riemannian map is totally godesic. We also obtain a decomposition theorem for the source manifold of an anti-invariant Riemannian map.\\

Let $F:(M_1,g_1) \longrightarrow (M_2,g_2)$ be a map and let
$p_2=F(p_1)$ for each $p_1 \in M_1$. Suppose that $\nabla^2$ is the
Levi-Civita connection on $(M_2,g_2)$. For $X \in \Gamma(TM_1)$ and
$V \in \Gamma(TM_2)$, we have
$${\nabla^F}^2_X(V\circ F)=\nabla^2_{F_*X}V,$$ where ${\nabla^F}^2$ is the pulback connection of $\nabla^2$. From now on, for simplicity, we
denote by $\nabla^2$ both the Levi-Civita connection of $(M_2, g_2)$
and its pullback along $F$. Then according to \cite{Nore}, for any
vector field $X$ on $M_1$ and any section $V$ of $(range
F_*)^\perp$, where $(range F_*)^\perp$ is the
 subbundle of $F^{-1}(TM_2)$ with fiber $(F_*(T_{p_1}M_1))^\perp$-orthogonal complement of $F_*(T_{p_1}M_1)$ for $g_2$ over ${p_1}$, we have
 $\nabla^{^{F \perp}}_XV$ which is the orthogonal projection of $\nabla^2_XV$ on $(F_*(TM_1))^\perp$. In \cite{Nore}, the author also showed that $\nabla^{^{F \perp}}$ is a linear connection on $(F_*(TM_1))^\perp$ such that $\nabla^{^{F \perp}}g_2=0$. We now suppose that $F$ is a Riemannian map and define $A_{V}$ as
\begin{equation}
\nabla^2_{_{F_*X}}V=-A_{_V}F_*X+\nabla^{^{F \perp}}_{_{X}}V,
\label{eq:4.1}
\end{equation}
where $A_{_V}F_*X$ is the tangential component (a vector field along
$F$) of $\nabla^2_{_{F_*X}}V$. Observe that $\nabla^2_{_{F_*X}}V$ is
obtained from the pullback connection of $\nabla^2$.  Thus, at
$p_1\in M_1$, we have $\nabla^2_{_{F_*X}}V(p_1) \in T_{F(p_1)}M_2$,
$A_{_V}F_*X(p_1) \in F_{*p_1}(T_{p_1}M_1)$ and $\nabla^{^{F
\perp}}_{_{X}}V(p_1)\in (F_{*p_1}(T_{p_1}M_1))^\perp$. It is easy to
see that $A_V F_*X$ is bilinear in $V$ and $F_*X$ and $A_V F_*X$ at
${p_1}$ depends only on $V_{p_1}$ and $F_{*{p_1}}X_{p_1}$. By direct
computations, we obtain
\begin{equation}
g_2(A_{_V} F_*X,F_*Y)=g_2(V, (\nabla F_*)(X,Y)), \label{eq:4.2}
\end{equation}
for $X, Y \in \Gamma((ker F_*)^\perp)$ and $V \in \Gamma((range
F_*)^\perp)$. Since $(\nabla F_*)$ is symmetric, it follows that
$A_{_V}$ is a symmetric linear transformation of $range F_*$.\\

Let $F:(M_1,g_1) \longrightarrow (M_2,g_2)$ be a map between
Riemannian manifolds $(M_1,g_1)$ and $(M_2,g_2)$, let $p_2=F(p_1)$
for each $p_1 \in M_1$.  Then the adjoint map ${^*F}_*$ of $F_*$ is
characterized by $g_1(x,{^*F}_{*p_1}y)=g_2(F_{*p_1}x,y)$ for $x \in
T_{p_1}M_1$, $y \in T_{F(p_1)}M_2$ and $p_1 \in M_1$. Considering
$F^{h}_*$ at each $p_1 \in M_1$ as a linear transformation
$$F^{h}_{*{p_1}}: ((ker F_*)^{\perp}(p_1),{g_1}_{{p_1}((ker F_*)^\perp (p_1))})
\rightarrow (range F_*(p_2),{g_2}_{{p_2}(range F_*)(p_2))}),$$ we
will denote the adjoint of $F^{h}_* $ by ${^*F^{h}}_{*p_1}.$ Let
${^*F}_{*p_1}$ be the adjoint of $F_{*p_1}:(T_{p_1}M_1,{g_1}_{p_1})
\longrightarrow (T_{p_2}M_2, {g_2}_{p_2})$. Then the linear
transformation
$$({^*F}_{*p_1})^h:range F_*(p_2) \longrightarrow (ker F_*)^{\perp}(p_1)$$
defined by $({^*F}_{*p_1})^hy= {^*F}_{*p_1}y$, where $y \in
\Gamma(range F_{*p_1}), p_2=F(p_1)$, is an isomorphism and
$(F^h_{*p_1})^{-1}=({^*F}_{*p_1})^h={^*(F^h_{*p_1})}$.\\

We are now ready to give necessary and sufficient conditions for
an anti-invariant Riemannian map from an almost Hermitian manifold
to a Riemannian manifold to be a totally geodesic map. We recall
that a differentiable map $F$ between Riemannian manifolds
$(M_1,g_1)$ and $(M_2,g_2)$ is called a totally geodesic map if
$(\nabla F_*)(X,Y)=0$ for all $X, Y \in \Gamma(TM_1)$.\\

\noindent{\bf Theorem 4.1.~}{\it Let $F$ be an anti-invariant
Riemannian map from a K\"{a}hler manifold $(M_1,g_1,J_1)$
to a Riemannian manifold $(M_2,g_2)$. Then $F$ is totally geodesic
if and only if}
\begin{equation}
{^*F}_{*}(A_VF_*(J_1X))\in \Gamma(\mu) \label{eq:4.3}
\end{equation}
\begin{equation}
{^*F}_{*}(A_VF_*(Z_1))\in \Gamma(J_1(ker F_*)) \label{eq:4.4}
\end{equation}
\begin{eqnarray}
g_2(F_*(J_1Y),(\nabla F_*)(X,\mathcal{B}Z))&=&g_2((\nabla
F_*)(X,J_1Y),F_*(\mathcal{C}Z))\label{eq:4.5}\\
g_1(\nabla_X\mathcal{B}Z,\mathcal{B}\bar{Z})&=&g_2((\nabla F_*)(X
,\mathcal{B}Z)+(\nabla F_*)(X ,\mathcal{C}Z),
F_*(\mathcal{C}\bar{Z}))\nonumber\\
&&-g_2(F_*(\mathcal{C}Z),(\nabla
F_*)(X,\mathcal{B}\bar{Z}))\label{eq:4.6}
\end{eqnarray}
{\it for $X,Y \in \Gamma(ker F_*)$, $Z_1 \in \Gamma(\mu)$, $Z,
\bar{Z}\in \Gamma((ker F_*)^\perp)$ and $V \in \Gamma((range
F_*)^\perp)$.}

\noindent{\bf Proof.~}From (\ref{eq:2.8}), (\ref{eq:2.9}) and
(\ref{eq:2.10}) we get
$$g_2((\nabla
F_*)(X,Z),F_*(\bar{Z}))=-g_1(\nabla_XJ_1Z,J_1\bar{Z})$$ for $X \in
\Gamma(ker F_*)$ and $Z, \bar{Z} \in \Gamma((ker F_*)^\perp)$.
Using (\ref{eq:3.1}) we have
\begin{eqnarray}
g_2((\nabla
F_*)(X,Z),F_*(\bar{Z}))&=&-g_1(\nabla_X\mathcal{B}Z,\mathcal{B}\bar{Z})-g_1(\nabla_X\mathcal{B}Z,\mathcal{C}\bar{Z})\nonumber\\
&-&g_1(\nabla_X\mathcal{C}Z,\mathcal{B}\bar{Z})-g_1(\nabla_X\mathcal{C}Z,\mathcal{C}\bar{Z}).\nonumber
\end{eqnarray}
Since $F$ is a Riemannian map, we obtain
\begin{eqnarray}
g_2((\nabla
F_*)(X,Z),F_*(\bar{Z}))&=&-g_1(\nabla_X\mathcal{B}Z,\mathcal{B}\bar{Z})-g_2(F_*(\nabla_X\mathcal{B}Z),F_*(\mathcal{C}\bar{Z}))\nonumber\\
&+&g_1(\mathcal{C}Z,\nabla_X\mathcal{B}\bar{Z})-g_2(F_*(\nabla_X\mathcal{C}Z),F_*(\mathcal{C}\bar{Z})).\nonumber
\end{eqnarray}
Then Riemannian map $F$ and (\ref{eq:2.8}) imply that
\begin{eqnarray}
g_2((\nabla
F_*)(X,Z),F_*(\bar{Z}))&=&-g_1(\nabla_X\mathcal{B}Z,\mathcal{B}\bar{Z})+g_2((\nabla F_*)(X,\mathcal{B}Z),F_*(\mathcal{C}\bar{Z}))\nonumber\\
&-&g_2(F_*(\mathcal{C}Z),(\nabla
F_*)(X,\mathcal{B}\bar{Z}))\nonumber\\
&+&g_2((\nabla
F_*)(X,\mathcal{C}Z),F_*(\mathcal{C}\bar{Z})).\label{eq:4.7}
\end{eqnarray}
In a similar way, one can obtain
\begin{eqnarray}
g_2((\nabla F_*)(X,Y),F_*(Z))&=&-g_2(F_*(J_1Y),(\nabla
F_*)(X,\mathcal{B}Z))\nonumber\\
&+&g_2((\nabla F_*)(X,J_1Y),F_*(\mathcal{C}Z)).\label{eq:4.8}
\end{eqnarray}
for $X,Y \in \Gamma(ker F_*)$ and $Z \in \Gamma((ker F_*)^\perp)$.
On the other hand, from (\ref{eq:2.8}) we have
$$g_2((\nabla F_*)(Z_1,Z_2),V)=g_2(\nabla^{^F}_{Z_1}F_*(Z_2),V).$$
for $Z_1,Z_2\in \Gamma(\mu)$ and $V\in \Gamma((range F_*)^\perp)$.
Hence, we get
$$g_2((\nabla F_*)(Z_1,Z_2),V)=-g_2(F_*(Z_2),\nabla^{^2}_{F_*(Z_1)}V).$$
Then, using (\ref{eq:4.1}) we obtain
$$g_2((\nabla F_*)(Z_1,Z_2),V)=g_2(F_*(Z_2),A_{V}F_*(Z_1)).$$
Thus we have
\begin{equation}
g_2((\nabla
F_*)(Z_1,Z_2),V)=g_1(Z_2,{^*F}_{*}(A_{V}F_*(Z_1))).\label{eq:4.9}
\end{equation}
In a similar way, we obtain
\begin{equation}
g_2((\nabla
F_*)(J_1X,J_1Y),V)=g_1(J_1Y,{^*F}_{*}(A_{V}F_*(J_1X)))\label{eq:4.10}
\end{equation}
for $X, Y\in \Gamma(ker F_*)$ and $V \in \Gamma((range
F_*)^\perp)$. Thus proof comes from
(\ref{eq:4.7})-(\ref{eq:4.10}).\\

Let $F$ be a Riemannian map from a Riemannian manifold to a
Riemannian manifold $(M_2,g_2)$. We say that a Riemannian map is a
Riemannian map with totally umbilical fibers if
\begin{equation}
h_2(X, Y)=g_1(X,Y)H \label{eq:4.11}
\end{equation}
for $X, Y \in \Gamma(ker F_*)$, where $h_2$ and  $H$ are the second fundamental form and  the mean curvature
vector field of the distribution $ker F_*$, respectively. In the sequel we show that anti-invariant Riemannian map puts some restrictions on the geometry of the distribution $ker F_*$.\\

\noindent{\bf Theorem~4.2~} {\it Let $F$ be a Lagrangian
Riemannian map with totally umbilical fibers  from a K\"{a}hler
manifold $(M_1,g_1,J)$ to a Riemannian manifold $(M_2,g_
)$. If $dim(ker F_{*p})>1$, $p\in M_1$, then the distribution $ker F_*$ defines a totally geodesic foliation on $M_1$.}\\

\noindent{\bf Proof.~} Since $F$ is a Lagrangian map with totally
umbilical fibres, for $X, Y\in \Gamma(ker F_*)$ and $Z \in
\Gamma((ker F_*)^\perp)$, we have
$g_1(\nabla^{^1}_XY,Z)=g_1(X,Y)g_1(H,Z)$. Using (\ref{eq:2.9}) we
obtain $g_1(\nabla^{^1}_XJY,JZ)=g_1(X,Y)g_1(H,Z)$. Hence we get
$$-g_1(\nabla^{^1}_XJZ,JY)=g_1(X,Y)g_1(H,Z).$$
Then (\ref{eq:4.11}) implies that
$$-g_1(X,JZ)g_1(H,JY)=g_1(X,Y)g_1(H,Z).$$
Thus we have
$$g_1(H,JY)JX=g_1(X,Y)H.$$
Taking inner product both sides with $JX$ and using (\ref{eq:2.9}), we arrive at
$$g_1(H,JY)g_1(X,X)=g_1(X,Y)g_1(H,JX).$$
Since $dim(ker F_*)>1
$, we can choose unit vector fields $X$ and $Y$ such that $g_1(X,Y)=0$, thus we derive
$$g_1(H,JY)g_1(X,X)=0.$$
Since $F$ is Lagrangian, the above equation implies that $H=0$ which shows that $ker F_*$ is totally geodesic.\\

Let us recall the notion of twisted product. Let $(M^{^m}, g_{_M})$ and $(N^{^n}, g_{_N})$ be Riemannian manifolds of dimensions $m$ and $n$, respectively. Let $P_1:M\times N \longrightarrow M$ and $P_2:M\times N \longrightarrow N$ be the canonical projections. Suppose that $f:M \times N \longrightarrow (0, \infty)$ be a smooth function. Then the twisted product of $(M^{^m}, g_{_M})$ and $(N^{^n}, g_{_N})$ with twisting function $f$ is defined to be the product manifold $\bar{M}=M \times N$ with metric tensor $\bar{g}=g_{_M} \oplus f^2 g_{_N}$ given by
\begin{equation}
\bar{g}=P^*_1g_{_M} + f^2 P^*_2 g_{_N}.\nonumber
\end{equation}
We denote this twisted product manifold $(\bar{M}, \bar{g})$ by $M \times_f N$. We now  recall that we
have the following result of \cite{Reckziegel}, {\sl Let $D_1$ be a vector subbundle
in the tangent bundle of a Riemannian manifold $M$ and $D_2$ be its
normal bundle. Suppose that the two distributions are involutive. We
denote the integral manifolds of $D_1$ and $D_2$ by $M_1$ and $M_2$,
respectively. Then $M$ is locally isometric to twisted product $M_1
\times_f M_2$ if the integral manifold $M_1$ is totally geodesic and
the integral manifold $M_2$ is totally umbilical, i.e, $M_2$  is a
totally
umbilical submanifold. }\\

Thus from the above remark and Theorem 4.2, we obtain the following result.\\

\noindent{\bf Corollary~4.1.~}{\it Let $M_1$ be a K\"{a}hler manifold and $M_2$ a Riemannian manifold. Then there do not exist a Lagrangian Riemannian map $F$ from $M_1$ to $M_2$ such that $M_1$ is a locally twisted product manifold of the form $M_{(ker F_*)^\perp}\times_f M_{(ker F_*)}$, where $M_{(ker F_*)^\perp}$ and $M_{(ker F_*)}$ are the integral manifolds of $(ker F_*)^\perp$ and $ker F_*$, and $f$ is the twisting function.}\\

Let $M_1$ be a K\"{a}hler manifold with complex structure $J$ and
$M_2$ a Riemannian manifold. We recall that a smooth map
$F:M_1\longrightarrow M_2$ is  called pluriharmonic if the second
fundamental form $\nabla F_*$ of the map $F$ satisfies $(\nabla
F_*)(X,Y)+(\nabla F_*)(JX,JY)= 0$ for any $X,Y\in \Gamma(TM_1)$,
\cite{Ohnita}. It is well known that if $M_1$ and $M_2$ are
K\"{a}hler manifolds
and $F:M_1\longrightarrow M_2$ is a holomorphic map, then $F$ is pluriharmonic. Since every pluriharmonic map is harmonic map, a holomorphic map is a harmonic map between K\"{a}hler manifolds.  For Lagrangian Riemannian maps from K\"{a}hler manifolds to Riemannian manifolds, we have the following result.\\

\noindent{\bf Theorem~4.3.~}{\it Let $(M_1,g_1,J)$ be a connected K\"{a}hler manifold and $M_2$ a Riemannian manifold. If a Langragian Riemannian map from $M_1$ to $M_2$ is pluriharmonic, then it is totally geodesic.}\\

\noindent{\bf Proof.} Suppose that $F$ is a pluriharmonic
Lagrangian Riemannian map from $M_1$ to $M_2$. Then we have
$$(\nabla F_*)(X,Y)+(\nabla F_*)(JX,JY)= 0$$ for  $X,Y\in \Gamma(ker F_*)$. On the other hand, from \cite{Sahin2}, we know that $$(\nabla F_*)(JX,JY)\in \Gamma((range F_*)^\perp),\quad (\nabla F_*)(X,Y) \in \Gamma(range F_*).$$
Since $TM_2=range F_* \oplus (range F_*)^\perp$, we get
\begin{equation}
(\nabla F_*)(JX,JY)=0 \label{eq:4.12}
\end{equation}
and
\begin{equation}
(\nabla F_*)(X,Y)=0. \label{eq:4.13}
\end{equation}
We also claim that $(\nabla F_*)(X,JY)=0$ for $X, Y \in \Gamma(ker F_*)$. Suppose that $F$ is pluriharmonic and $(\nabla F_*)(X,JY)\neq 0$. But since $F$ is pluriharmonic, we have $(\nabla F_*)(X,JY)-(\nabla F_*)(JX,Y)=0.$
Then using (\ref{eq:2.8}) we have
$$-F_*(\nabla^{^1}_XJY)+F_*(\nabla^{^1}_YJX)=0.$$
Since $M_1$ is K\"{a}hler manifold, we derive
$$-F_*(J\nabla^{^1}_XY)+F_*(J\nabla^{^1}_YX)=0.$$
Hence we obtain
$$F_*(J[X,Y])=0.$$
Hence, we conclude that $J[X,Y]\in \Gamma(ker F_*)$ which implies that $[X,Y] \in \Gamma((ker F_*)^\perp)$. This tells that $ker F_*$ is not integrable which contradicts with the rank theorem. Since $F$ is a subimmersion, it follows that the rank of $F$ is constant on $M_1$, then the rank theorem for functions implies that $ker F_*$ is an integrable subbundle of $TM_1$, (\cite{Marsden}, page:205). Thus we should have
\begin{equation}
(\nabla F_*)(X,JY)=0\label{eq:4.14}
\end{equation}
for $X, Y \in \Gamma(ker F_*)$. Then proof comes from (\ref{eq:4.12})-(\ref{eq:4.14}).\\


\begin{thebibliography}{25}
\bibitem{Marsden} Abraham,R., Marsden, J.E., Ratiu,T., {\it Manifolds, Tensor Analysis, and Applications}, Springer-Verlag, Newyork, 1988.
\bibitem{Baird-Wood} Baird, P. and  Wood, J. C., {\it Harmonic
Morphisms Between Riemannian Manifolds}, London Mathematical
Society Monographs, No. 29, Oxford University Press, The Clarendon
Press, Oxford, 2003.

\bibitem{Domingo} Chinea, D., Almost contact metric
submersions. Rend. Circ. Mat. Palermo, (1985), 34(1), 89-104.


\bibitem{Rham} De Rham, G., Sur la reducibilite dun espace de Riemann. Comment.
Math. Helv. 26,(1952), 328–344.


\bibitem{Escobales} Escobales, R. H. Jr., Riemannian
submersions from complex projective space. J. Differential Geom.,
(1978), 13(1), 93-107.

\bibitem{Falcitelli-Ianus-Pastore} Falcitelli, M.,
Ianus, S., Pastore, A. M., {\it Riemannian Submersions and Related
Topics}. World Scientific, River Edge, NJ, 2004.

\bibitem{Fischer} Fischer, A. E.: Riemannian maps between
Riemannian manifolds, Contemporary math. 132, 331-366, (1992).

\bibitem{Garcia-Rio-Kupeli}Garcia-Rio, E., Kupeli,D. N.,  {\it Semi-Riemannian maps and their
Applications}, Kluwer Academic, Dortrecht, 1999.

\bibitem{Gray} Gray, A., Pseudo-Riemannian almost product
manifolds and submersions, J. Math. Mech., (1967), 16, 715-737.

\bibitem{Ianus-Mazzocco-Vilcu} Ianus, S., Mazzocco, R.,
Vilcu, G. E., Riemannian submersions from quaternionic manifolds.
Acta Appl. Math., (2008), 104(1), 83-89.

\bibitem{Marrero-Rocha} Marrero, J. C., Rocha, J., Locally
conformal K\"{a}hler submersions. Geom. Dedicata, (1994), 52(3),
271-289.

\bibitem{Nore}Nore, T., Second fundamental form of a map, Ann. Mat. Pur. and Appl., 146, (1987),  281-310.

\bibitem{Ohnita}Ohnita, Y., On pluriharmonicity of stable harmonic maps, J. London Math. Soc. (2) 35 (1987) 563-568.

 \bibitem{O'Neill} O'Neill, B., The fundamental
equations of a submersion, Mich. Math. J., (1966), 13, 458-469.
\bibitem{Reckziegel} Ponge R., Reckziegel H., Twisted products in pseudo-Riemannian geometry, Geom. Dedicata, 1993, 48(1), 15–25.

\bibitem{Sahin1}\d{S}ahin, B., Anti-invariant
Riemannian submersions from almost Hermitian manifolds, Central
European J.Math, 8(3), (2010), 437-447.

\bibitem{Sahin2}\d{S}ahin, B., Invariant and anti-invariant
Riemannian maps to K\"{a}hler manifolds, International Journal of
Geometric Methods in Modern Physics, vol:7, no:3 (2010), 1-19.

\bibitem{Watson} Watson, B., Almost Hermitian submersions. J.
Differential Geometry, (1976), 11(1), 147-165.

\bibitem{Yano-Kon} Yano, K. and Kon, M., {\it Structures on
Manifolds}, World Scientific, Singapore, 1984.

\end{thebibliography}
\end{document}